\documentclass[titlepage,11pt]{article}
\usepackage{latexsym}
\usepackage{amsfonts}
\usepackage{amsmath}
\newcommand\blackslug{\hbox{\hskip 1pt \vrule width 4pt height 8pt depth 1.5pt
        \hskip 1pt}}
\newcommand\bbox{\hfill \quad \blackslug \bigbreak}

\def\LL{,\ldots,}

\title{A shorter proof of the path-width theorem}
\author{
Paul Seymour\thanks{Supported by AFOSR grant
FA9550-22-1-0234, and NSF grant  DMS-2154169.}\\
Princeton University, Princeton, NJ 08544}

\date{September 9, 2023; revised \today}

\newtheorem{thm}{}[section]

\newcommand{\Proof}{\noindent{\bf Proof.}\ \ }

\begin{document}
\maketitle
\begin{abstract}
A graph has {\em path-width} at most $w$ if it can be built from a sequence of graphs each with at most $w+1$ vertices,
by overlapping consecutive terms. Every graph with path-width at least $w-1$ 
contains every $w$-vertex forest as a minor: this was originally proved by Bienstock, Robertson, Thomas and the author,
and was given a short proof by Diestel. Here we give a proof even shorter and simpler than that of Diestel.
\end{abstract}

\section{The proof}
All graphs in this paper are finite, and may have loops or parallel edges. If $G$ is a graph, $|G|$ denotes its number of vertices, and for $A\subseteq V(G)$, $G[A]$ denotes the subgraph induced on $A$.
A {\em path-decomposition} of a graph $G$ is a sequence $(W_1\LL W_n)$ of subsets of $V(G)$ (called {\em bags}), with union $V(G)$,
such that for every edge $uv$ of $G$ there exists $i$ such that $u,v\in W_i$, and such that $W_i\cap W_k\subseteq W_j$
for $1\le i<j<k\le n$; and it has {\em width} at most $w$ if $|W_i|\le w+1$ for each $i$.  
A graph has {\em path-width} at most $w$ if it admits a path-decomposition with width at most $w$. 
Robertson and the author~\cite{GM1}
proved that for every forest $F$, all graphs that do not contain $F$ as a minor have bounded path-width 
(and the conclusion is false for all graphs $F$ that are not forests); and later Bienstock, Robertson, Thomas and author~\cite{bienstock} proved:
\begin{thm}\label{mainthm}
For every forest $F$, every graph that does not contain $F$ as a minor has path-width at most $|F|-2$.
\end{thm}
This is tight, since a complete graph on $|F|-1$ vertices has path-width $|F|-2$ and does not contain
$F$ as a minor. It was given a short proof by Diestel~\cite{diestel}, but there is an even shorter proof, that 
we present here.

A {\em model} of a loopless graph $H$ in a graph $G$ is a map $\phi$ with domain $V(H)\cup E(H)$, such that
\begin{itemize}
\item $\phi(h)$ is a non-null connected subgraph of $G$ for each $h\in V(H)$, and $\phi(h),\phi(h')$ are vertex-disjoint
for all distinct $h,h'\in V(H)$;
\item $\phi(f)\in E(G)$ for each $f\in E(H)$, and $\phi(f)\ne \phi(f')$ for all distinct $f,f'\in E(H)$;
\item if $f\in E(H)$ is incident in $H$ with $h\in V(H)$, then $\phi(f)$ is incident in $G$ with a vertex of $\phi(h)$.
\end{itemize}
Thus there is a model of $H$ in $G$ if and only if $G$ contains $H$ as a minor.

A {\em separation} of $G$ is a pair 
$(A,B)$
of subsets of $V(G)$ with union $V(G)$, such that there are no edges between $A\setminus B$ and $B\setminus A$, and 
its {\em order} is $|A\cap B|$. If $(A,B)$ and $(A',B')$ are separations of $G$, we write $(A,B)\le (A',B')$ if 
$A\subseteq A'$ and $B'\subseteq B$. For each integer $w\ge 0$, we say a separation $(A,B)$ of a graph $G$ is 
{\em $w$-good}
if there is a path-decomposition of $G[A]$ with width at most $w$ and with last bag $A\cap B$.
We need the following observation, which is the heart of the proof:
\begin{thm}\label{extend}
If $(A',B')$ and $(P,Q)$ are separations of $G$, where $(A',B')$ is $w$-good and $(P,Q)\le (A',B')$, and 
there are $|P\cap Q|$ vertex-disjoint paths of $G$ between $P$ and $B'$, then $(P,Q)$ is $w$-good.
\end{thm}
\Proof
Let $t=|P\cap Q|$, and let $R_1\LL R_t$ be disjoint paths between $P$ and $B'$. We may assume that each has only one 
vertex in $B'$, and hence in $A'\cap B'$.  Each of these paths has only its first vertex in $P$, and so 
if we contract the edges of $R_1\LL R_t$, we preserve the subgraph $G[P]$. 
Let $H$ be the union of $G[P]$ and the paths $R_1\LL R_t$. Since $(A',B')$ is $w$-good, there is a path-decomposition of 
$H$ of width at most $w$, such that its last bag consists of the $t$ ends in $B'$ of the paths $R_1\LL R_t$. 
But contracting
the edges of $R_1\LL R_t$ brings this to a path-decomposition of $G[P]$ with last bag $P\cap Q$ (since each edge to 
be contracted has both ends inside a bag). This proves \ref{extend}.~\bbox

If $(A,B)$ and $(A',B')$ are separations of $G$, the second {\em extends} the first if $(A,B)\le (A',B')$  
and $|A\cap B|\ge |A'\cap B'|$.
A $w$-good separation of $G$ is {\em maximal} if no different $w$-good separation extends it.
Let $w\ge 0$ be an integer, let $T$ be a tree or the null graph, and let $(A,B)$ be a separation of a graph $G$. We say that
$(A,B)$ is {\em $(w,T)$-spanning} if
\begin{itemize}
\item $|A\cap B|=|T|$;
\item there is a model $\phi$ of $T$ in $G[A]$ such that $V(\phi(h))\cap A\cap B\ne \emptyset$ for each $h\in V(T)$; and
\item if $|T|\le w+1$ then $(A,B)$ is maximal $w$-good.
\end{itemize}
In order to prove \ref{mainthm}, we may assume that $F$ is a tree $T$ say (by adding edges to $F$ if necessary), and so 
it suffices to prove:
\begin{thm}\label{excellentthm}
Let $w\ge 0$ be an integer, let $G$ be a graph that has path-width more than $w$, and let $T$ be a tree or the null graph,
with $|T|\le w+2$.
Then there is a $(w,T)$-spanning separation of $G$.
\end{thm}
\Proof
We proceed by induction on $|T|$, keeping $w$ fixed. If $|T|=0$, the result holds since 
there is a maximal $w$-good separation of order zero, say $(A,B)$ (possibly with $A=\emptyset$), which is therefore
$(w,T)$-spanning. So we assume that $1\le |T|\le w+2$ and the result holds for $|T|-1$. Choose $j\in V(T)$ with degree at most one,
and if $|T|\ge 2$ let $i$ be the neighbour of $j$ in $T$. 

From the
inductive hypothesis, there is a $(w, T\setminus \{j\})$-spanning separation $(A,B)$ of $G$, which is therefore 
maximal $w$-good, since $|T\setminus \{j\}|<w+2$. 
Let $\phi$ be a model of
$T\setminus \{j\}$ in $G[A]$ such that $V(\phi(h))\cap A\cap B\ne \emptyset$ for each $h\in V(T)\setminus \{j\}$.
We choose $v\in B\setminus A$ as follows.
If $|T|=1$, then $A\cap B=\emptyset$; choose $v\in B$ arbitrarily. (This is possible since $B\ne \emptyset$, because $G$ has path-width more 
than $w$: this is the only place where we use that the path-width is large.) If $|T|\ge 2$, 
let $u\in V(\phi(i))\cap  B$. Then $u$ has a neighbour $v\in B\setminus A$, since otherwise $(A,B\setminus \{u\})$ is $w$-good
and extends $(A,B)$, contradicting the maximality of $(A,B)$. This defines $v$. 

If $|T|=w+2$, then $(A\cup \{v\},B)$ is $(w,T)$-spanning, so we may assume that $|T|<w+2$,
and therefore $(A\cup \{v\},B)$ is $w$-good.
So there
is a maximal $w$-good separation $(A',B')$ of $G$ that extends $(A\cup \{v\},B)$. Since $(A',B')$ does not extend $(A,B)$
(because $(A,B)$ is maximal $w$-good), its order is exactly $|T|$.
Suppose that there is a separation $(P,Q)$ of $G$ of order less than $|T|$, with $(A\cup \{v\},B)\le (P,Q)\le (A',B')$.
Choose $(P,Q)$ with minimum order; then it follows from Menger's theorem that there are $|P\cap Q|$ vertex-disjoint
paths from $P$ to $B'$, and so from \ref{extend}, $(P,Q)$ is $w$-good. But $(P,Q)$ extends $(A,B)$, since 
$|P\cap Q|\le |T|-1= |A\cap B|$, and $(P,Q)\ne (A,B)$ since $v\in P$, contradicting
the maximality of $(A,B)$.
Thus there is no such $(P,Q)$, and so by Menger's theorem, there are $|T|$ disjoint paths of $G$ between $A\cup \{v\}$ 
and $B'$. By combining these with the model $\phi$, we deduce that $(A',B')$ is $(w,T)$-spanning. This proves \ref{excellentthm}.~\bbox

\end{document}